\documentclass[12pt]{amsart}
\usepackage{amsmath,amscd,amssymb,amsfonts,graphics,mathrsfs}
\usepackage{hyperref}
\newcommand{\hl}{\hyperlink}
\newcommand{\htt}{\hypertarget}
\newcommand{\h}{\hbox}
\newcommand{\q}{\quad}

\newcommand{\bs}{\par\bigskip}
\newcommand{\ms}{\par\medskip}
\newcommand{\sk}{\par\smallskip}
\newcommand{\bsn}{\par\bigskip\noindent}
\newcommand{\msn}{\par\medskip\noindent}
\newcommand{\skn}{\par\smallskip\noindent}
\newcommand{\ges}{\geqslant}
\newcommand{\les}{\leqslant}
\newcommand{\1}{\hskip1pt}

\newcommand{\mcup}{\hbox{$\bigcup$}}

\newcommand{\D}{{\mathscr D}}

\newcommand{\F}{{\mathscr F}}
\newcommand{\Fb}{\F^{\ssb}}
\newcommand{\Hc}{{\mathscr H}}
\newcommand{\I}{{\mathscr I}}
\newcommand{\J}{{\mathscr J}}
\newcommand{\Lc}{{\mathscr L}}

\newcommand{\M}{{\mathscr M}}
\newcommand{\Sc}{{\mathscr S}}
\newcommand{\OO}{{\mathscr O}}
\newcommand{\X}{{\mathscr X}}
\newcommand{\Y}{{\mathscr Y}}
\newcommand{\A}{{\mathbb A}}
\newcommand{\PP}{{\mathbb P}}
\newcommand{\Q}{{\mathbb Q}}
\newcommand{\C}{{\mathbb C}}

\newcommand{\RR}{{\mathbf R}}
\newcommand{\DD}{{\mathbf D}}
\newcommand{\Z}{{\mathbb Z}}
\newcommand{\kk}{{\mathbf k}}

\newcommand{\ep}{\varepsilon}

\newcommand{\Ht}{\widetilde{H}}

\newcommand{\Sb}{{}\,\overline{\!S}{}}

\newcommand{\DR}{{\rm DR}}
\newcommand{\Gr}{{\rm Gr}}
\newcommand{\lcd}{{\rm lcd}}
\newcommand{\dd}{\partial}
\newcommand{\pHc}{{}^{\bf p}\!\Hc}
\newcommand{\RHDQ}{{\rm RHD}_{\Q}}
\newcommand{\Hom}{{\mathscr H}\hskip-1.7pt om}
\newcommand{\Ext}{{\mathscr E}\hskip-1.4pt xt}
\newcommand{\bl}{\bigl}
\newcommand{\br}{\bigr}
\newcommand{\pl}{\1{+}\1}
\newcommand{\mi}{\1{-}\1}
\newcommand{\eq}{\,{=}\,}
\newcommand{\caps}{\,{\cap}\,}
\newcommand{\ins}{\,{\in}\,}
\newcommand{\tos}{\,{\to}\,}
\newcommand{\sst}{\,{\subset}\,}
\newcommand{\stm}{\,{\setminus}\,}
\newcommand{\less}{\,{\leqslant}\,}
\newcommand{\gess}{\,{\geqslant}\,}

\newcommand{\ssb}{\raise.15ex\h{${\scriptscriptstyle\bullet}$}}
\newcommand{\ssc}{\,\raise.15ex\h{${\scriptstyle\circ}$}\,}

\newcommand{\into}{\hookrightarrow}
\newcommand{\simto}{\,\,\rlap{\hskip1.5mm\raise1.4mm\hbox{$\sim$}}\hbox{$\longrightarrow$}\,\,}
\newcommand{\indlim}{\rlap{\raise-6.2pt\h{$\,\to$}}{\rm lim}}
\begin{document}
\title[Local cohomological dimension]{Topological calculation of local
  cohomological dimension}
\author[T.~Reichelt]{Thomas Reichelt}
\address{Thomas Reichelt : Lehrstuhl f\"ur Algebraische Geometrie,
  Universit\"at Mannheim, B6 26, 68159 Mannheim, Germany}
\email{reichelt@math.uni-mannheim.de}
\author[M.~Saito]{Morihiko Saito}
\address{Morihiko Saito : RIMS Kyoto University, Kyoto 606-8502, Japan}
\email{msaito@kurims.kyoto-u.ac.jp}
\author[U.~Walther]{Uli Walther}
\address{Uli Walther : Purdue University, Dept.\ of Mathematics, 150 N.\ University St., West Lafayette, IN 47907, USA}
\email{walther@math.purdue.edu}
\thanks{TR was supported by DFG Emmy-Noether-Fellowship RE 3567/1-2. MS was supported by JSPS Kakenhi 15K04816. UW was supported by NSF grant DMS-2100288 and by Simons Foundation Collaboration Grant for Mathematicians \#580839.}
\begin{abstract} We show that the sum of the local cohomological dimension and the rectified $\mathbb Q$-homological depth of a closed analytic subspace of a complex manifold coincide with the dimension of the ambient manifold. The local cohomological dimension is then calculated using the cohomology of the links of the analytic space. In the algebraic case the first assertion is equivalent to the coincidence of the rectified $\mathbb Q$-homological depth with the de~Rham depth studied by Ogus, and follows essentially from his work. As a corollary we show that the local cohomological dimension of a quasi-projective variety is determined by that of its general hyperplane section together with the link cohomology at 0-dimensional strata of a complex analytic Whitney stratification.
\end{abstract}
. 
\maketitle
\centerline{\bf Introduction}
\bsn
Let $Y$ be a reduced closed subvariety of a smooth complex algebraic variety $X$. The {\it local cohomological dimension\1} $\lcd(Y,X)$ of $Y\sst X$ (\cite{Ogu}) is defined as
\htt{1}{}
$$\max\bl\{j\ins\Z\mid\Hc^j_Y(\F)\,{\ne}\,0\,\,\h{for some quasi-coherent}\,\,\F\br\}.
\leqno(1)$$
One can show that this maximum does not change by restricting $\F$ to the structure sheaf $\OO_X$ (using the inductive limit of coherent subsheaves on affine open subvarieties together with free resolutions of these coherent subsheaves and also the spectral sequence associated with the truncations $\sigma_{\ges\ssb}$, see \hl{S1.4}{1.4} below). This notion comes from a problem raised by A.~Grothendieck (see \cite[p.\,79]{Ha0}) where he notes that it can be bounded by the number of equations describing $Y\sst X$ locally (using the \v Cech cohomology). We have $\lcd(Y,X)\gess{\rm codim}_XY$, since the equality holds if $Y$ is smooth (or more generally a local complete intersection).
\sk
The above definition is extended to the analytic case where $Y$ is a closed analytic subspace of a complex manifold $X$, using the {\it algebraic local cohomology\1} $\Hc^j_{[Y]}$ (see \cite{Gr2}, \cite{Me}, \cite{Ka2}, and also \hl{S1.3}{1.3} below), that is,
$$\lcd(Y,X):=\max\bl\{j\ins\Z\mid\Hc^j_{[Y]}(\OO_X)\ne0\br\}.$$
From now on, we assume that $X,Y$ are {\it analytic.}
\sk
We can calculate $\lcd(Y,X)$ in a {\it purely topological way\1} using the $t$-structure on $D^b_c(X,\C)$ constructed in \cite{BBD}. It is not very difficult to prove the following (see \hl{S2.1}{2.1} below and also \cite{RSW}).
\par\htt{T1}{}\msn
{\bf Theorem~1.} {\it In the above notation, we have the equality
\htt{2}{}
$$\aligned\lcd(Y,X)={}&d_X-\RHDQ(Y)\q\q\h{with}\\ \RHDQ(Y):={}&\min\bl\{j\ins\Z\mid{}\pHc^j(\Q_Y)\ne0\br\},\endaligned
\leqno(2)$$
in particular, $d_X{-}\,\lcd(Y,X)$ depends only on $Y$ as in the algebraic case.}
\sk
Here $d_X\,{:=}\,\dim X/\C$, ${}\pHc^{\ssb}$ is the cohomological functor associated with the self-dual $t$-structure for the middle perversity on $D^b_c(X,\C)$ constructed in \cite{BBD}, and $\RHDQ(Y)$ is called the {\it rectified $\Q$-homological depth,} see \cite[1.1.3]{HL} and also \cite[p.\,387]{Sch}. The notion of {\it rectified homotopical depth,} denoted by ${\rm rhd}(Y)$, was introduced in \cite{Gr2}. As is remarked above for $\lcd(Y,X)$, we have
\htt{3}{}
$$\RHDQ(Y)\les d_Y,
\leqno(3)$$
looking at smooth points of $Y$. (Indeed, ${}\pHc^j\Q_Y\eq 0$ for $j\,{\ne}\,d_Y$ if $Y$ is smooth.)
\sk
Employing Thom's isotopy lemma (see for instance \cite[Prop.\,11.1]{Ma}), Theorem~\hl{T1}{1} implies the following.
\par\htt{C1}{}\msn
{\bf Corollary~1.} {\it The local cohomological dimension $\lcd(Y,X)$ is invariant under a topologically equisingular deformation of $(Y,X)$ assuming $Y$ is compact or has only isolated singularities.}
\ms
Here a {\it topologically equisingular deformation\1} of $(Y,X)$ means that there is a smooth morphism $p\,{:}\,\X\tos T$ (with $T$ smooth) together with a complex analytic Whitney stratification of a closed analytic subspace $\Y\sst\X$ whose strata are smooth over $T$ with $(\Y_0,\X_0)\eq(Y,X)$ for some $0\ins T$. We assume $p|_{\Y}\,{:}\,\Y\tos T$ is proper if $Y$ is compact. The following is already known.
\par\htt{T2}{}\msn
{\bf Theorem~2} (Hamm, L\^e \cite[Thm.\,1.9]{HL}). {\it For a complex analytic space $Y$ and $k\ins\Z_{>0}$, the following conditions are equivalent\,$:$
\skn
{\rm\,(i)\1} $\RHDQ(Y)\ges k$.
\skn
{\rm(ii)} For any locally closed irreducible analytic subset $Z\subset Y$, we have}
\htt{4}{}
$$\Hc^j_Z\Q_Y\eq 0\q\q(\forall\,j<k\mi{\dim Z}).
\leqno(4)$$
\ms
The last condition should be compared with (\hl{1.1.2}{1.1.2}) below (where we may assume that the Whitney stratification is compatible with $Z$ locally on $Y$) and also with the following definition of the {\it de~Rham depth\1} ${\rm DRD}(Y)$ of $Y$ in the algebraic setting (see \cite[2.12]{Ogu}):
\msn
{\bf Definition 1.} We have ${\rm DRD}(Y)\ges k$ $\iff$ For any (not necessarily closed) point $y\in Y$,
\htt{5}{}
$$H^j_y(Y)\eq 0\q\q(\forall\,j<k\mi{\dim\overline{\{y\}}}).
\leqno(5)$$
\skn
Indeed, in the algebraic case, it is shown (see \cite[Thm.\,2.13]{Ogu}) that
\htt{6}{}
$$\lcd(Y,X)=\dim X-{\rm DRD}(Y).
\leqno(6)$$
Theorem~\hl{T1}{1} is then equivalent to the following.
\par\htt{C2}{}\msn
{\bf Corollary~2.} {\it For a complex algebraic variety $Y$, we have the equality}
\htt{7}{}
$${\rm DRD}(Y)=\RHDQ(Y^{\rm an}).
\leqno(7)$$
\ms
The definition of $H^j_y(Y)$ in \cite{Ha2} employs the {\it completion\1} along $Y$ of the de~Rham complex of a smooth ambient variety $X$. Consequently the translation to our setting is not trivial. However, the comparison with singular cohomology has been shown in \cite[Ch.\,4,\,Thm.\,1.1]{Ha2} when $y$ is a {\it closed\1} point. Moreover the restriction to the closed points is allowed by {\it constructibility\1} of the subset $S_iY$ in \cite[2.15.1]{Ogu}. (These reasonings do not seem sufficiently clear to the reader in the paper.) Note that the latter subset corresponds to $S_i^{\vee}(\Q_Y)$ in (\hl{1.1.7}{1.1.7}) below. Theorem~\hl{T1}{1} in the algebraic case then follows essentially from \cite[Thm.\,2.13 combined with 2.15.2]{Ogu}. It is indeed stated on p.\,328 in that paper as follows:
\sk
``It is reasonable to expect that, more generally, the singularities of $Y$ determine lcd$(X, Y)$ whenever $X$ is smooth, and indeed this is the case. Our main result expresses lcd$(X, Y)$ in terms of the classically studied cohomology groups $H^i(Y, Y - y, \C)=H^i_y(Y, \C)$, studied for instance by Milnor [1968]." (Unfortunately it does not seem to be mentioned explicitly what ``our main result" means.)
\sk
In the complex analytic case, Theorem~\hl{T1}{1} gives the following.
\par\htt{C3}{}\msn
{\bf Corollary~3.} {\it For a closed analytic subspace $Y$ in a complex manifold $X$ and $k\ins\Z$, the following conditions are equivalent\,$:$}
\skn
\rlap{\,\rm(i)}\hskip.7cm\hangindent=.7cm
$\lcd(Y,X)\les d_X\mi k$,\,\,\,{\it that is,}\,\,\,$\RHDQ(Y)\ges k$.
\skn
\rlap{\rm(ii)}\hskip.7cm\hangindent=.7cm
{\it The reduced cohomology $\Ht^j\bl(L(T_S,y_S),\Q\br)$ of the link $L(T_S,y_S)$ vanishes for $j\less k\mi d_S\mi 2$ and for any stratum $S$ of a complex analytic Whitney stratification $\Sc$ of $\,Y$ with $d_S\,{<}\,k$, where $T_S$ is a local transversal slice to $S$ with $\{y_S\}\eq S\caps T_S$.}
\ms
We assume all the strata of $\Sc$ are connected. The {\it link\1} $L(Z,z)$ at $z\ins Z$ of a complex analytic space $Z$ in general is defined to be $Z\caps S_z$ with $S_z$ a sufficiently small sphere in an ambient smooth space with center $z$. (The well-definedness follows from the argument in the proof of the topological cone theorem using integral curves of a controlled vector field.) Corollary~\hl{C3}{3} follows from Theorem~\hl{T1}{1} together with local topological triviality along strata of a complex analytic Whitney stratification \cite[Cor.\,10.6]{Ma}. In the case $T_S\eq\{y_S\}$ and $L(T_S,y_S)\eq\emptyset$, we have $\Ht^j(\emptyset,\Q)\eq\Q$ for $j\eq{-}1$, and $0$ otherwise. (This is compatible with the isomorphism (\hl{2.2.1}{2.2.1}) below.) Condition~(ii) then implies that any irreducible component of $Y$ has dimension $\ges k$ (setting $j\eq{-}1$ in the inequality $j\less k\mi d_S\mi 2$). By Corollary~\hl{C3}{3} for $k\eq0$, we always have
\htt{8}{}
$$\RHDQ(Y)\ges 0.
\leqno(8)$$
In the algebraic case, the inequality $\lcd(Y,X)\less d_X$ follows from Grothendieck's vanishing theorem, see for instance \cite[6.1.2]{BS}. Observe that the equality $\lcd(Y,X)\eq d_X$ holds if there is an {\it isolated point\1} of $Y$. Corollary~\hl{C3}{3} for $k\eq1$ gives the following.
\par\htt{C4}{}\msn
{\bf Corollary~4.} {\it For $X,Y$ as in Corollary~{\rm\hl{C3}{3}}, we have $\lcd(Y,X)\less d_X\mi 1$, that is, $\RHDQ(Y)\gess 1$, if and only if $\,Y$ does not contain an isolated point.}
\sk
This corresponds to a special case of a theorem of Hartshorne-Lichtenbaum (see \cite[Thm.\,3.1]{Ha1}, where smoothness of $X$ is not required), see also \cite[Cor.\,2.10]{Ogu}, \cite[Cor.\,1.19]{MP}.
\sk
Assuming $Y$ has no isolated points, we have the equality $\RHDQ(Y)\eq 1$ if $Y$ contains a 1-dimensional irreducible component or if $Y\eq Y_1\cup Y_2$ with $Y_1,Y_2$ smooth, $d_{Y_1},d_{Y_2}\gess 2$, and $Y':=Y_1\cap Y_2$ is non-empty and 0-dimensional. For the last case, we can use the long exact sequence associated with the cohomological functor ${}\pHc^{\ssb}$ applied to the distinguished triangle
$$\Q_Y\to\Q_{Y_1}{\oplus}\,\Q_{Y_2}\to\Q_{Y'{}}\buildrel{+1}\over\to.$$
The assertion in the last case follows also from Corollary~\hl{C5}{5} just below (which is a special case of Corollary~\hl{C3}{3} with $k\eq2$).
\par\htt{C5}{}\msn
{\bf Corollary~5.} {\it For $X,Y$ as in Corollary~{\rm\hl{C3}{3}}, the following conditions are equivalent\,$:$}
\skn
\rlap{\,\rm(i)}\hskip.7cm\hangindent=.7cm
$\lcd(Y,X)\les d_X\mi 2$,\,\,\,{\it that is,}\,\,\,$\RHDQ(Y)\ges 2$,
\skn
\rlap{\rm(ii)}\hskip.7cm\hangindent=.7cm
{\it Any irreducible component of $Y$ has dimension at least $2$, and the link $L(Y,y)$ is connected for any $y\in Y$.}
\ms
Indeed, using the topological cone theorem, the link $L(Y,y)$ in (ii) can be replaced by $Y\caps B_y\stm\{y\}$ (with $B_y$ a sufficiently small ball in an ambient smooth space). The latter is connected if it is non-empty and $y$ is in a positive-dimensional stratum (using local topological triviality along the stratum). The assertion corresponding to Corollary~\hl{C5}{5} has been studied quite well in the algebraic case, see \cite[Thm.\,7.5]{Ha1} for cones of projective varieties, \cite[Cor.\,5.5]{PS} for the positive characteristic case where Frobenius depth is used instead of de~Rham depth, \cite[Cor.\,2.11]{Ogu} for the characteristic zero case, \cite[Thm.\,2.9]{HL2} for varieties over a separably closed field where invariance by completion \cite[Prop.\,2.2]{Ha1} is used, and \cite[Thm.\,1.4]{Zh} for unramified regular local rings of mixed characteristic. The assertion for varieties over a field $\kk$ does not necessarily hold unless $\kk$ is separably closed (or {\it strict\1} Henselization is used when one considers the connectivity of the complement of the closed point in ${\rm Spec}\,R$ with $R$ a local ring), see also Example~\hl{E1.4}{1.4} below. For extensions to the case $\lcd(Y,X)\less d_X\mi 3$, see for instance \cite{Va}, \cite{DT} (and also \cite{WZ}).
\sk
Returning to our situation, Corollary~\hl{C3}{3} gives also the following.
\par\htt{C6}{}\msn
{\bf Corollary~6.} {\it For a closed subvariety $Y$ of a quasi-projective smooth complex variety $X$ and $k\ins\Z_{>0}$, the following conditions are equivalent\,$:$}
\skn
\rlap{\,\rm(i)}\hskip.7cm\hangindent=.7cm
$\lcd(Y,X)\les d_X\mi k$,\,\,\,{\it that is,}\,\,\,$\RHDQ(Y)\ges k$.
\skn
\rlap{\rm(ii)}\hskip.7cm\hangindent=.7cm
{\it $\lcd(Y',X')\less d_X\mi k$, that is, $\RHDQ(Y')\gess k\mi 1$, for a general hyperplane section $X'\sst X$ with $Y'\,{:=}\,Y\caps X'$, and $\Ht^j\bl(L(Y,y),\Q\br)\eq 0$ for $j\less k\mi 2$ and for any point $y$ in a $0$-dimensional stratum of a complex analytic Whitney stratification of $Y$.}
\ms
Here $X,Y$ are viewed as analytic spaces, and a general hyperplane section means that the hyperplane section $X'$ is smooth and transversal to any stratum of a complex analytic Whitney stratification of $Y$ (hence $X'$ does not meet any 0-dimensional stratum). Corollary~\hl{C6}{6} says that the local cohomological dimension of a quasi-projective variety can be calculated by using a general hyperplane section provided that one has information about this invariant (topological or algebraic) {\it at worst singular points.} The maximum in (\hl{1}{1}) is usually attained at the worst singular points. In the cone case we have the following.
\par\htt{C7}{}\msn
{\bf Corollary~7.} {\it Assume $X\eq\C^{\,d_X}$, and $Y$ is the affine cone of a complex projective variety $Z\sst\PP^{\,d_X-1}$. Then for $k\ins\Z_{>0}$, the following conditions are equivalent\,$:$}
\skn
\rlap{\,\rm(i)}\hskip.7cm\hangindent=.7cm
{\it $\lcd(Y,X)\less d_X{-}k$,\,\,\,that is,\,\,\,$\RHDQ(Y)\gess k$.}
\skn
\rlap{\rm(ii)}\hskip.7cm\hangindent=.7cm
{\it $\lcd(Z,\PP^{\,d_X-1})\less d_X{-}k$, that is, $\RHDQ(Z)\gess k{-}1$, and $H^j(Z)/H^j(\PP^{\,d_X-1})\eq0$ for $j\less k{-}2$ with $\Q$-coefficients.}
\ms
In the case $Z$ is smooth and connected, condition~(ii) is equivalent to the condition that $H^j_{\rm prim}(Z)=0$ for $1\less j\less k{-}2$ and $k\less d_Y$. If $Z$ is the projective cone of a smooth connected projective variety $V\sst\PP^{d_X-2}$, we see that the maximum in (\hl{1}{1}) is attained at the generator of the cone which is the pull-back of the vertex of the cone $Z$ by the rational map from $Y$ to $Z$. (Indeed, $\Ht^j(Z)\eq H^{j-2}(V)(-1)$ for any $j\ins\Z$, see (\hl{2.4.3}{2.4.3}) below.)
\sk
In Section~\hl{S1}{1} we review some basics of $t$-structure, holonomic $\D$-modules, algebraic local cohomology, and local cohomological dimension. In Section~\hl{S2}{2} we prove Theorem~\hl{T1}{1} and Corollaries~\hl{C3}{3} and \hl{C6}{6}--\hl{C7}{7} using the assertions in Section~\hl{S1}{1}.
\bs\bs\htt{S1}{}
\vbox{\centerline{\bf 1. Preliminaries}
\bsn
In this section we review some basics of $t$-structure, holonomic $\D$-modules, algebraic local cohomology, and local cohomological dimension.}
\par\htt{S1.1}{}\msn
{\bf 1.1.~$t$-structure.} Let $X$ be a reduced complex analytic space. We denote by $D^b_c(X,\Q)$ the derived category of bounded $\Q$-complexes on $X$ with constructible cohomology sheaves. For $k\ins\Z$ we have the full subcategories $D^b_c(X,\Q)^{\les k}$, $D^b_c(X,\Q)^{\ges k}$ defined by the following conditions for $\Fb\in D^b_c(X,\Q)$ respectively\,:
\htt{1.1.1}{}
$$\Hc^j\Fb|_S=0\q\q(j>k\mi d_S),
\leqno(1.1.1)$$
\vskip-6mm
\htt{1.1.2}{}
$$\,\,\,\,\Hc^j_S\Fb=0\q\q(j<k\mi d_S).
\leqno(1.1.2)$$
Here $S$ runs over strata of a complex analytic Whitney stratification $\Sc$ of $X$ compatible with $\Fb$, and $d_S:=\dim S/\C$, see \cite[2.1.2]{BBD} (where $p(S)\eq{-}d_S$) and \cite{Di}, \cite{KS}, \cite{Sch}.
This definition is independent of the Whitney stratification (as is seen just below).
Recall that
\htt{1.1.3}{}
$$\Hc^j_S\Fb=\Hc^ji_S^!\,\Fb\q\q(j\ins\Z),
\leqno(1.1.3)$$
with $i_S:S\into X$ the inclusion, see (\hl{1.1.10}{1.1.10}) below. Notice that a shift of indices by $d_S$ occurs in the last conditions of (\hl{1.1.1}{1.1.1}--\hl{1.1.2}{2}) if we restrict $\Fb$ to a local transversal slice to each stratum $S\ins\Sc$ as in the proof of Corollary~\hl{C3}{3}, see \hl{S2.2}{2.2} below. There is also a shift by $2d_S$ (or $-2d_S$) if we apply $\RR\Gamma_{\{x\}}$ (or the dual functor $\DD$) to a local system on $S$ with $x\in S$. (This is closely related to \cite[2.14]{Ogu}.)
\sk
We then see that the above two conditions are respectively equivalent to the following {\it support\1} and {\it cosupport\1} conditions (implying the independence of stratification):
\htt{1.1.4}{}
$$\dim\,\overline{S_i(\Fb)}\les k-i\q(\forall\,i\ins\Z),
\leqno(1.1.4)$$
\vskip-6mm
\htt{1.1.5}{}
$$\dim\,\overline{S^{\vee}_i(\Fb)}\les i-k\q(\forall\,i\ins\Z),
\leqno(1.1.5)$$
with
\htt{1.1.6}{}
$$S_i(\Fb):=\bl\{x\in X\mid\Hc^i\Fb_x\ne 0\br\},
\leqno(1.1.6)$$
\vskip-5mm
\htt{1.1.7}{}
$$\,S^{\vee}_i(\Fb):=\bl\{x\in X\mid\Hc^i_{\{x\}}\Fb\ne 0\br\},
\leqno(1.1.7)$$
see \cite{Di}, \cite{KS}, \cite{Sch}. (The cosupport condition (\hl{1.1.5}{1.1.5}) seems to be related closely to \cite[2.14--15]{Ogu}.)
\sk
For $k\ins\Z$, set
$$D^b_c(X,\Q)^{[k]}:=D^b_c(X,\Q)^{\les k}\cap D^b_c(X,\Q)^{\ges k}.$$
These are abelian full categories of $D^b_c(X,\Q)$, see \cite{BBD}. By definition this $t$-structure is {\it self-dual,} that is,
\htt{1.1.8}{}
$$\DD\bl(D^b_c(X,\Q)^{\les k}\br)=D^b_c(X,\Q)^{\ges-k},
\leqno(1.1.8)$$
in particular, $D^b_c(X,\Q)^{[0]}$ is stable by the functor $\DD$ assigning the dual.
Recall that, in the $X$ {\it smooth\1} case, we have (omitting the Tate twist)
\htt{1.1.9}{}
$$\DD(\Fb)=\RR\Hom\1_{\Q}\bl(\Fb,\Q_X[2d_X]\br).
\leqno(1.1.9)$$
\par\htt{R1.1a}{}\msn
{\bf Remark~1.1a.} The above argument is valid with $\Q$ replaced by any subfield of $\C$.
\par\htt{R1.1b}{}\msn
{\bf Remark~1.1b.} For a closed analytic subset $Y\subset X$, we have the isomorphisms
\htt{1.1.10}{}
$$\RR\Gamma_Y=(i_Y)_*\1i_Y^!,\q i_Y^!=\DD\ssc i_Y^*\ssc\DD,
\leqno(1.1.10)$$
with $i_Y:Y\into X$ the canonical inclusion. Indeed, if $j':X\stm Y\into X$ denotes the inclusion, we have the distinguished triangles
\htt{1.1.11}{}
$$\aligned(i_Y)_*\1i_Y^!&\to{\rm id}\to\RR j'_*\1j'{}^{-1}\buildrel{+1}\over\to\,,\\ \RR\Gamma_Y&\to{\rm id}\to\RR j'_*\1j'{}^{-1}\buildrel{+1}\over\to\,.\endaligned
\leqno(1.1.11)$$
These imply the first isomorphism in (\hl{1.1.10}{1.1.10}) (non-canonically). The second isomorphism can be reduced (non-canonically) to
\htt{1.1.12}{}
$$j'_!=\DD\ssc\RR j'_*\ssc\DD,
\leqno(1.1.12)$$
since $\DD^2={\rm id}$. (Some more argument is required to get the {\it canonical\1} isomorphisms due to the {\it ambiguity of mapping cones.})
\par\htt{R1.1c}{}\msn
{\bf Remark~1.1c.} In the case $\Fb=\Q_Y$ with $Y$ a closed analytic subset of $X$, we have
\htt{1.1.13}{}
$$S^{\vee}_i(\Q_Y)=\bl\{y\in Y\mid\Ht^{i-1}(S_y\caps Y,\Q)\ne 0\br\},
\leqno(1.1.13)$$
where $S_y$ is a sufficiently small sphere in $X$ with center $y$ (using a topological cone theorem together with (\hl{1.1.11}{1.1.11})). This may be related to the work of Milnor \cite{Mi} in the sentence quoted from \cite{Ogu} in the introduction.
\par\htt{R1.1d}{}\msn
{\bf Remark~1.1d.} The {\it topological cone theorem\1} asserts a stratified homeomorphism
\htt{1.1.14}{}
$$Y\caps B_y\cong{\rm Cone}(Y\caps S_y).
\leqno(1.1.14)$$
Here $S_y:=\dd B_y$ with $B_y$ a sufficiently small closed ball in an ambient smooth space with center $y$, and Cone$(E)$ for a topological space $E$ means a topological cone, which is obtained by collapsing $E{\times}\{0\}$ in $E{\times}[0,1]$. Moreover the distance function $\delta_y$ from $y$ on the left-hand side coincides with the second projection on the right-hand side up to a constant multiple. (This is proved by using the integral curves of a controlled $C^{\infty}$ vector field $\xi$ on $B_y\stm\{y\}$ satisfying $\langle\1\xi,{\rm d}\delta_y\rangle=1$.)
\sk
The topological cone theorem is rather nontrivial. It is easier to show that for a bounded constructible complex $K^{\ssb}\in D^b_c(X,\Q)$, the cohomology groups $H^j(S_{x,\ep},K|_{S_{x,\ep}})$ ($j\ins\Z$) do not depend on the radius $\ep$ of sufficiently small spheres $S_{x,\ep}$ in an ambient complex manifold with center $x\in X$ fixed. We can use for instance {\it embedded resolutions\1} of $(X^{(k)},X^{(k-1)})$ with $X^{(k)}$ the union of strata of dimension $\les k$ of a complex analytic Whitney stratification of $X$. (Here it is easy to construct a controlled vector field using a partition of unity.) We can replace $S_{x,\ep}$ with $B_{x,\ep}\stm B^{\circ}_{x,\ep'}$ for $0\,{<}\,\ep'\,{<}\,\ep\,{\ll}\,1$, where $B^{\circ}_{x,\ep'}$ is the interior of $B_{x,\ep'}$.
\par\htt{R1.1e}{}\msn
{\bf Remark~1.1e.} For a closed analytic subset $Y\subset X$ and $k\ins\Z_{>0}$, we have by using (\hl{1.1.2}{1.1.2})
\htt{1.1.15}{}
$$\RHDQ(Y)\ges k\iff \Q_Y\in D^b_c(X,\Q)^{\ges k}.
\leqno(1.1.15)$$
\par\htt{S1.2}{}\msn
{\bf 1.2.~Holonomic $\D$-modules.} Let $X$ be a complex manifold, and $\D_X$ be the sheaf of holomorphic differential operators. It has the filtration $F$ by order of differential operators such that the associated graded ring $\Gr^F_{\ssb}\D_X$ is locally isomorphic to the polynomial ring over $\OO_X$, and ${\rm Specan}_X\Gr^F_{\ssb}\D_X$ is naturally identified with the cotangent bundle $T^*\!X$, that is, $\Gr^F_{\ssb}\D_X$ is identified with the direct image to $X$ of the sheaf of holomorphic functions on $T^*\!X$ whose restriction to each fiber is a polynomial.
\sk
A coherent left $\D_X$-module $\M$ is called {\it holonomic\1} if there is locally an increasing filtration $F$ on $\M$ such that $(\M,F)$ is a filtered $\D_X$-module, $F_p\1\M=0$ for $p\ll0$, $\Gr^F_{\ssb}\M$ is a coherent $\Gr^F_{\ssb}\D_X$-module, and ${\rm Supp}\,\Gr^F_{\ssb}\M\subset T^*\!X$ has dimension $d_X:=\dim X$.
\sk
We denote by $M_{\rm hol}(\D_X)$ the category of holonomic left $\D_X$-modules. This is an abelian subcategory of the category of left $\D_X$-modules $M(\D_X)$, which is closed by subquotients and extensions in $M(\D_X)$. Moreover it is stable by the dual functor $\DD$.
\par\htt{R1.2a}{}\msn
{\bf Remark~1.2a.} For a holonomic $\D_X$-module $\M$, let $\DR_X(\M)$ be the de~Rham complex shifted by $d_X$. It is well known that
\htt{1.2.1}{}
$$\DR_X(\M)\in D^b_c(X,\C)^{[0]}.
\leqno(1.2.1)$$
Indeed, it is shown in \cite{Ka1} that
\htt{1.2.2}{}
$$\DR_X(\M)\in D^b_c(X,\C)^{\les 0}.
\leqno(1.2.2)$$
The assertion (\hl{1.2.1}{1.2.1}) follows from (\hl{1.2.2}{1.2.2}) using (\hl{1.1.8}{1.1.8}), since it is well known that the de~Rham functor $\DR_X$ is compatible with the dual functor $\DD$ (see \cite[Prop.\,2.4.12]{mhp} for the filtered case).
\sk
As a corollary of (\hl{1.2.1}{1.2.1}), we get the functorial isomorphisms
\htt{1.2.3}{}
$$\DR_X\ssc\Hc^j={}\pHc^j\ssc\DR_X\q\q(j\ins\Z),
\leqno(1.2.3)$$
where ${}\pHc^{\ssb}$ is the cohomological functor associated with the $t$-structure on $D^b_c(X,\C)$ which is constructed in \cite{BBD}.
\par\htt{R1.2b}{}\msn
{\bf Remark~1.2b.} In the case $\M$ is a {\it regular holonomic\1} $\D_X$-module, it is also possible to prove (\hl{1.2.1}{1.2.1}) by induction on $d\eq{\dim{\rm Supp}\,\M}$ using \cite{De}. Here the regular holonomic $\D$-modules can be defined by induction on $\dim{\rm Supp}\,\M$ using \cite{De}, see \cite{rh}.
\par\htt{S1.3}{}\msn
{\bf 1.3.~Algebraic local cohomology.} Let $X$ be a complex manifold with $Y\subset X$ a closed analytic subset. The {\it algebraic local cohomology\1} $\Hc^j_{[Y]}\M$ for an $\OO_X$-module $\M$ is defined by
\htt{1.3.1}{}
$$\Hc^j_{[Y]}\M:=\rlap{\raise-10.5pt\h{$\,\,\,\scriptstyle k$}}\indlim\,\Ext^j_{\OO_X}(\OO_X/\I_Y^k,\M)\q\q(j\ins\Z),
\leqno(1.3.1)$$
with $\I_Y\subset\OO_X$ the ideal of $Y$, see \cite{Gr2}.
This can be extended to the case of a $\D_X$-module $\M$ using an injective resolution of $\M$ as $\D_X$-modules, since
\htt{1.3.2}{}
$$\Hom_{\OO_X}(\OO_X/\I_Y^k,\M)=\M^{\I_Y^k},
\leqno(1.3.2)$$
with right-hand side the subsheaf of $\M$ annihilated by $\I_Y^k$ (and $\dd_{x_i}\I_Y^k\subset\I_Y^{k-1}$ with $x_i$ local coordinates), see also \cite{Ka2}. Recall that an injective $\D_X$-module is an injective $\OO_X$-module, since $\D_X$ is flat over $\OO_X$.
\sk
It is known that the $\Hc^j_{[Y]}\M$ are {\it holonomic\1} if so is $\M$ (where {\it regularity\1} is not assumed), see \cite{Ka2}. Here $b$-functions in a generalized sense are used.
\par\htt{R1.3a}{}\msn
{\bf Remark~1.3a.} We can define also $\RR\Gamma_{[Y]}$ by replacing $\Ext^j_{\OO_X}(\OO_X/\I_Y^k,\M)$ in (\hl{1.3.1}{1.3.1}) with
$$\Hom_{\OO_X}(\OO_X/\I_Y^k,\J^{\ssb}),$$
where $\J^{\ssb}$ is an injective resolution of $\M$. It is well known (see for instance \cite{Me}) that
\htt{1.3.3}{}
$$\DR_X\RR\Gamma_{[Y]}\OO_X=\RR\Gamma_Y\DR_X\OO_X,
\leqno(1.3.3)$$
as a consequence of \cite[Lem.\,17]{AH} (applied to an embedded resolution), see also \cite[Thm.\,2]{Gr1}. This is closely related to {\it regularity\1} (although the equivalence of categories is not needed for our purpose).
\par\htt{R1.3b}{}\msn
{\bf Remark~1.3b.} For the calculation of the local cohomology, we have locally a \v Cech-type complex associated with local generators $g_j$ ($j\ins J$) of the ideal of $Y$ and defined by using the $\Hc^0_{[X|Z_{\!J'}]}$ ($J'\sst J$). Here $Z_{\!J'}\,{:=}\,\1\mcup_{j\in J'}\1Z_j$ with $Z_j\sst X$ the hypersurface defined by $g_j$, and $\Hc^0_{[X|Z_{\!J'}]}$ is the localization along $Z_{\!J'}$, see \cite[Thm.\,A1.3]{Eis} for stalks. This can be applied to the proof of the commutativity of the local cohomology with the direct image of $\D$-modules by an embedded resolution of $Y$, see \cite{rh}. Using this argument, we can prove also the {\it holonomicity of the local cohomology\1} in the structure sheaf case by reducing to the normal crossing divisor case.
\par\htt{S1.4}{}\msn
{\bf 1.4.~Local cohomological dimension.} For a closed subvariety $Y$ of a smooth variety $X$ over a field $\kk$, the local cohomological dimension (\cite{Ogu}) is defined by
\htt{1.4.1}{}
$$\lcd(X,Y):=\max\bl\{j\ins\Z\mid\Hc^j_Y(\F)\ne0\,\,\,(\exists\,\F\,\,\h{quasicoherent)}\br\}.
\leqno(1.4.1)$$
Here it is enough to consider the case $\F=\OO_X$. Indeed, the assertion is essentially local, and we may assume $X={\rm Spec}\,A$, where a quasicoherent sheaf $\F$ corresponds to an $A$-module $M$. Since the local cohomology commutes with inductive limit, and any $A$-module is the inductive limit (or union) of its finite $A$-submodules, it is enough to consider the case of coherent sheaves.
\sk
Let $\Lc^{\ssb}\to\F$ be a free resolution, where the $\Lc^p$ are free sheaves of finite length, and vanish for $p>0$. We have the spectral sequence (associated with the truncations $\sigma_{\ges\ssb}$)
\htt{1.4.2}{}
$$E_1^{p,q}=\Hc^q_Y\Lc^p\Longrightarrow\Hc^{p+q}_Y\F.
\leqno(1.4.2)$$
This implies that it is sufficient to consider the structure sheaf case. The above argument may be called in \cite{Ogu} ``d\'evissage" (which is sometimes used in the Grothendieck school).
\par\htt{E1.4}{}\msn
{\bf Example~1.4} (with $\lcd(Y,X)\eq d_X{-}1$). Assume $\kk\eq\Q$, $X\eq\A^4_{\kk}$, and $Y_{\C}\,{:=}\,Y{\otimes}_{\kk}\C\sst\C^4$ is the union of $Y'_{\C}$ and its Galois (or complex) conjugate $Y''_{\C}$, where $Y'_{\C}\sst\C^4$ is defined by the equations
\htt{1.4.3}{}
$$i\1x^2\pl y^2\pl z^2=y^2\mi z^2\pl i\1w^2=0.
\leqno(1.4.3)$$
This is defined over the subfield $\kk'\,{:=}\,\kk[i]\sst\C$ with $i\eq\sqrt{-1}$. The Galois group ${\rm Gal}(\kk'/\kk)$ is isomorphic to ${\mathfrak S}_2\eq\{\pm 1\}$. (Here one may use Weil's theory with $\C$ the universal domain.) We see that $Y$ is irreducible over $\kk$, but $Y'_{\C}\cap Y''_{\C}=\{0\}$, hence $\Hc^{\1d_X-1}_Y\OO_X\ne 0$ using the base change by $\kk\into\C$. Note that $Y\sst X$ is an irreducible component of a complete intersection, but not a local complete intersection. 
\bs\bs\htt{S2}{}
\vbox{\centerline{\bf 2. Proof of the main theorem and corollaries}
\bsn
In this section we prove Theorem~\hl{T1}{1} and Corollaries~\hl{C3}{3} and \hl{C6}{6}--\hl{C7}{7} using the assertions in Section~\hl{S1}{1}.}
\par\htt{S2.1}{}\msn
{\bf 2.1.~Proof of Theorem~\hl{T1}{1}.} There are isomorphisms in $D^b_c(X,\C)$\,:
\htt{2.1.1}{}
$$\aligned&\DR_X\RR\Gamma_{[Y]}\OO_X=\RR\Gamma_Y\DR_X\OO_X=\RR\Gamma_Y\C_{X}[d_X]\\={}&i_Y^!\C_{X}[d_X]=\DD\1i_Y^*\1\DD(\C_{X}[d_X])=(\DD\1\C_Y)[-d_X],\endaligned
\leqno(2.1.1)$$
where the zero extension $(i_Y)_*$ is omitted to simplify the notation. Indeed, the first, the third and fourth, and the last isomorphisms follow from (\hl{1.3.3}{1.3.3}), (\hl{1.1.10}{1.1.10}), and (\hl{1.1.9}{1.1.9}) respectively. Since the de~Rham functor $\DR_X$ is faithful on holonomic $\D_X$-modules, the assertion then follows from (\hl{1.1.8}{1.1.8}) and (\hl{1.2.3}{1.2.3}), see also \cite[Section~1.1]{RSW}. This finishes the proof of Theorem~\hl{T1}{1}.
\par\htt{S2.2}{}\msn
{\bf 2.2.~Proof of Corollary~\hl{C3}{3}.} We have {\it local topological triviality\1} along each stratum $S$, see \cite[Cor.\,10.6]{Ma}. This implies that $Y$ is locally homeomorphic to $T_S{\times}U_{y_S}$ as a stratified space. Here $T_S$ denotes a local {\it transversal slice\1} to $S$, and $U_{y_S}\sst S$ is a small open ball around $\{y_S\}\eq T_S\caps S$. Note that $\Q_Y$ is locally isomorphic to the pull-back of $\Q_{T_S}$ by the smooth projection $T_S{\times}U_{y_S}\to T_S$. Under the pull-back by this smooth projection with relative {\it real\1} dimension $2d_S$, the indices in the last conditions of (\hl{1.1.1}{1.1.1}--\hl{1.1.2}{2}) are {\it shifted by\1} $d_S$, since $p(S)$ in \cite[2.1.2]{BBD} is given by the real dimension of $S$ divided by 2 up to sign, see \cite[p.\,131]{Di}.
\sk
The assertion then follows from Theorem~\hl{T1}{1} using condition~(\hl{1.1.2}{1.1.2}), since we have the isomorphisms
\htt{2.2.1}{}
$$\Ht^{j-1}\bl(L(T_S,y_S),\Q\br)=\Hc^j_{\{y_S\}}\Q_{T_S}\q(\forall\,j\ins\Z),
\leqno(2.2.1)$$
using the topological cone theorem. Indeed, (\hl{2.2.1}{2.2.1}) holds even in the case $T_S\eq\{y_S\}$, and the reduced cohomology in (ii) vanishes for $j\less k\mi d_S\mi 2$ if $d_S\gess k$. This finishes the proof of Corollary~\hl{C3}{3}.
\par\htt{S2.3}{}\msn
{\bf 2.3.~Proof of Corollary~\hl{C6}{6}.} We can verify that a general hyperplane section intersects {\it any\1} positive-dimensional stratum $S$. (Indeed, it intersects the closure of $S$, and the intersection $Y'\caps\Sb$ cannot be contained in the boundary of $S$ using the transversality for the strata $S'\subset\Sb$ with $d_{S'}\eq d_S\mi 1$, where we may assume $Y$ is projective.) The assertion then follows from Corollary~\hl{C3}{3}, since one can get a transversal slice by intersecting general hyperplanes. This finishes the proof of Corollary~\hl{C6}{6}.
\par\htt{S2.4}{}\msn
{\bf 2.4.~Proof of Corollary~\hl{C7}{7}.} By Corollary~\hl{C6}{6} it is enough to show that
\htt{2.4.1}{}
$$\Ht^j\bl(L(Y,0)\br)\eq0\,\,\h{for}\,\,j\less k{-}2\iff H^j(Z)/H^j(\PP^{\,d_X-1})\eq0\,\,\,\h{for}\,\,j\less k{-}2.
\leqno(2.4.1)$$
Here the cohomology groups are with $\Q$-coefficients. We may assume that $k\less d_Y$ by (\hl{3}{3}). The injectivity of the restriction morphism $H^j(\PP^{\,d_X-1})\tos H^j(Z)$ for $j\less 2d_Z$ can be shown taking a general linear subspace of $\PP^{\,d_X-1}$ (that is, an intersection of general hyperplanes) which intersects $Z$ transversally at smooth points, and looking at its pull-back to a desingularization of $Z$. (One can also use the integration of the pull-back of a power of a K\"ahler form on $\PP^{\,d_X-1}$.)
\sk
Set $Y^{\circ}:=Y\stm\{0\}$ (with $0$ the vertex of the cone $Y$). This is homotopy equivalent to the link $L(Y,0)$, and has a structure of $\C^*$-bundle over $Z$. We then get the long exact sequence
\htt{2.4.2}{}
$$\to H^{j-2}(Z)(-1)\buildrel{\eta}\over\to H^j(Z)\to H^j(Y^{\circ})\to H^{j-1}(Z)(-1)\to,
\leqno(2.4.2)$$
where $\eta$ is the first Chern class of the $\C^*$-bundle, see for instance \cite[\S1.3]{RSW}. 
Using this, the equivalence (\hl{2.4.1}{2.4.1}) can be verified by increasing induction on $k$. This finishes the proof of Corollary~\hl{C7}{7}.
\par\htt{R2.4}{}\msn
{\bf Remark~2.4.} In the case $Z$ is a projective cone of a projective variety $V\sst\PP^{d_X-2}$, the complement $Z^{\ssc}$ of the vertex $0$ in $Z$ is a line bundle over $V$. We then get that
\htt{2.4.3}{}
$$\Ht^j(Z)=H^j_c(Z^{\ssc})=H^{j-2}(V)(-1)\q(j\ins\Z).
\leqno(2.4.3)$$

\end{document}